\documentclass[10pt]{article}

\usepackage{amssymb}
\usepackage{amsbsy}
\usepackage{graphicx}
\usepackage{amsmath}

\usepackage{times}

\begin{document}

\title{For what number of cars must self organization occur in the Biham-Middleton-Levine traffic model
from any possible starting configuration?}
\author{Tim D. Austin  \and Itai Benjamini}
\date{}

\maketitle

\begin{abstract}
For any initial configuration of fewer than $\frac{1}{2}N$ cars
the BML model will self organize to attain speed one. On the other
hand, there is a configuration of size $m$ in which no car can
move if and only if $m$ is at least $2N$.
\end{abstract}


\newenvironment{nmath}{\begin{center}\begin{math}}{\end{math}\end{center}}

\newtheorem{thm}{Theorem}
\newtheorem{qu}[thm]{Question}
\newtheorem{lem}[thm]{Lemma}
\newtheorem{prop}[thm]{Proposition}
\newtheorem{cor}[thm]{Corollary}
\newtheorem{conj}[thm]{Conjecture}
\newtheorem{dfn}[thm]{Definition}
\newtheorem{prob}[thm]{Problem}


\newcommand{\A}{\mathcal{A}}
\newcommand{\B}{\mathcal{B}}
\renewcommand{\Pr}{\mathbb{P}}
\newcommand{\s}{\sigma}
\renewcommand{\P}{\mathcal{P}}
\renewcommand{\O}{\Omega}
\renewcommand{\o}{\omega}
\newcommand{\e}{\varepsilon}
\renewcommand{\S}{\Sigma}
\newcommand{\T}{\mathrm{T}}
\newcommand{\co}{\mathrm{co}}
\renewcommand{\i}{\mathrm{i}}
\renewcommand{\l}{\lambda}
\newcommand{\U}{\mathcal{U}}
\newcommand{\G}{\Gamma}
\newcommand{\g}{\gamma}
\renewcommand{\L}{\Lambda}
\newcommand{\hcf}{\mathrm{hcf}}
\newcommand{\F}{\mathcal{F}}
\renewcommand{\a}{\alpha}
\renewcommand{\b}{\beta}
\newcommand{\Zn}{\bb{Z}_n}

\newcommand{\bb}[1]{\mathbb{#1}}
\renewcommand{\rm}[1]{\mathrm{#1}}
\renewcommand{\cal}[1]{\mathcal{#1}}

\newcommand{\qed}{\nolinebreak \hspace{\stretch{1}} $\Box$}

Write $\bb{Z}_N = \bb{Z} / (N\bb{Z})$, and consider the $N \times
N$ discrete torus $\bb{Z}_N \times \bb{Z}_N$.  We think of the
first coordinate as pointing vertically and the second as pointing
horizontally. Suppose we choose $m < N^2$ points of this torus at
random, and then assign to each of them a colour -- red or blue --
independently with equal probabilities. We think of these points
as little cars on the torus.  Having once been positioned, the
cars move at discrete time-steps according to the following rules:
\begin{enumerate}
\item Red cars try to move one step upwards, and blue cars try to move one step
to the right.
\item However, two cars may not occupy the same space.  In any
given time step, first all the blue cars try to move, then all the
red cars.  If, when the time comes to move, the space immediately
to the right of a blue car is already occupied by a red car, then
the blue car -- and any more blue cars in a line to the left of it
-- cannot move, so stay still.  Then, if any vertical line of red
cars is blocked by a blue car above them, they all stay still.
\end{enumerate}

In the classical Biham-Middleton-Levine (BML) model (\cite{BML}),
the number $m$ is taken to be $pN^2$ for some $p \in (0,1)$.
Numerical simulations suggest that the system exhibits a phase
transition: there is a certain critical value of $p$, say
$p_{\rm{c}}$ (thought to be roughly $0.3$), such that for $p
> p_{\rm{c}}$, asymptotically almost surely every car is blocked
from moving by some other car after finitely many steps (we say
the system is \textbf{stuck}), whereas for $p < p_{\rm{c}}$,
asymptotically almost surely the system self-organizes: every car
is blocked from moving only finitely many times, and so after some
finite time the cars all move freely (we say the system attains
\textbf{speed one}).

This behaviour has been witnessed in a number of computer
simulations for increasingly large $N$.  However, recently a more
subtle picture has started to emerge: intermediate phases of
behaviour have been observed between attaining speed one and
getting stuck, and it also now appears possible that there is no
$p_{\rm{c}} \in (0,1)$, but rather that as $N \to \infty$ the
behaviour changes increasingly sharply around a point
$p_{\rm{c}}(N)$ that decays very slowly to $0$. See, in
particular, D'Souza \cite{D'S}. However, to date the only result
to be proved rigorously (by Angel, Holroyd and Martin \cite{AHM})
is that the system gets stuck for all $p$ sufficiently large.
Furthermore, when the implicit lower bound for such $p$ in this
last result is computed it is found to be very close to $1$.

In this work we consider a related deterministic question, which
turns out to be much simpler: Given $N$, for which $m$ does
\emph{any} initial configuration of the cars inevitably attain
speed one, and for which is there some initial configuration of
the cars which gets stuck?

\begin{prop}\label{prop:detBML}
If $m < \frac{1}{2}N$, the system must attain speed one.
\end{prop}

\textbf{Remark}\hspace{5pt} Of course, this bound is very much
lower than the threshold ($\sim pN^2$) originally observed for the
random model.

\textbf{Proof}\hspace{5pt} The idea here is to consider the
locations of the cars on the $N$ NW-SE diagonals of the torus. Let
these be $D_1,D_2,\ldots,D_N$:
\[D_k = \{(i,j)\in \bb{Z}_N^2:\ i + j = k \mod N\}.\]
If a car moves during a time step, then it moves up one diagonal;
else it stays still. For each $t \geq 0$ let $\phi_t:\bb{Z}_N^2
\to \bb{Z}_N$ be the associated `time-corrected diagonal map':
$\phi(i,j) = i+j-t \mod N$.

Suppose $X^1,X^2,\ldots,X^m$ are the initial positions of our
cars, and write $X^i_t$ for the position of car $i$ at time $t
\geq 0$ (so $X^i_0 = X^i$) and $Y^i_t = \phi_t(X^i_t)$. Thus,
knowing the $Y^i_t$ at a given time $t$ tells us something about
the configuration $(X^i_t)_{i \leq m}$, but far from specifies it
uniquely. Our proof will use constraints on the behaviour of
$Y^1_t,Y^2_t,\ldots,Y^m_t$ as $t$ increases.

At a given time $t$, the points $Y^i_t$ are distributed within
$\bb{Z}_N$: some points of $\bb{Z}_N$ may be occupied by many such
$Y^i_t$ (several cars may occupy different points on the same
diagonal), while others will be empty.  We will partition the set
$\bb{Z}_N \setminus \{Y^i_t:\ i \leq m\}$ of empty points at time
$t$ into a union of arcs in $\bb{Z}_N$, say $A^1_t \cup A^2_t \cup
\cdots \cup A^{r(t)}_t$, where (for sake of argument) we label the
arcs in order and choose $A^1_t$ to be the arc containing the
first non-occupied point when $\bb{Z}_N$ is written as
$\{0,1,\ldots,N\}$.

\begin{lem}\label{lem:arcstatic}
Suppose that at time $t$ some arc $A^s_t = \{y,y+1,\ldots,y+l\}$
has length at least $2$. Then at time $t+1$ the point immediately
to its left, $y-1$, is still occupied by some $Y^i_{t+1}$, and the
set $A^s_t \setminus \{y+l\} = \{y,y+1,\ldots,y+l-1\}$ is still an
arc of unoccupied points.
\end{lem}

\textbf{Proof of Lemma \ref{lem:arcstatic}}\hspace{5pt} This is a
direct observation from the dynamics of the cars.  The lower
boundary point cannot move, since given a number of cars all in
the same diagonal in the discrete torus with no cars in the
diagonal above, at least one of those cars will not be blocked
during the next time step; so, under the time-corrected diagonal
map, at least one of the images of those cars must stay still in
this time step. Also, images of cars under the maps $\phi$ can
only either stay still or move one step to the left in one time
step, so it is clear that the points $\{y,y+1,\ldots,y+l-1\}$
cannot become occupied during the next time step. \qed

\begin{lem}\label{lem:arclength}
The dynamics cannot create new arcs $A^s_t$ of length greater than
$1$: the number of such long arcs is non-increasing in $t$.
\end{lem}

\textbf{Proof of Lemma \ref{lem:arclength}}\hspace{5pt} Suppose
that at time $t+1$ we have an arc $A^s_{t+1}$ of length at least
$2$. During the time step from $t$ to $t+1$, the images of those
cars that are now in diagonals immediately above and below
$A^s_{t+1}$ either stayed still or moved one step to the left.
Thus, by Lemma \ref{lem:arcstatic}, at time $t$ there must have
been an empty arc at least as long as $A^s_{t+1}$ and with the
same lower end-point. Thus to each empty arc of length at least
$2$ at time $t+1$ we can associate such an arc at time $t$; since
this association is also clearly unique, the number of such arcs
cannot increase. \qed

\begin{lem}\label{lem:nolongarcs}
If the system never attains speed one, then there must come a time
when no arcs $A^s_t$ have length greater than $1$.
\end{lem}

\textbf{Proof of Lemma \ref{lem:nolongarcs}}\hspace{5pt} The point
is that if there are infinitely many times at which cars are
blocked, then in particular some car must be blocked infinitely
often.  Suppose it is car $i \leq m$.  This means that as $t$
increases $Y^i_t$ describes infinitely many circuits around the
discrete circle $\bb{Z}_N$.

But after completing one such circuit (say at time $T$) , there
can be no arcs of length greater than $1$ remaining; for if there
is still an arc $A^s_T$ of length at least $2$ at time $T$, then,
arguing as in the proof of the previous lemma, there must be a
sequence of arcs $A^{s(t)}_t$ for $t = 0,1,\ldots,T$, all of
length at least $2$, and all with the same lower end-point, say
$y$. Thus we deduce that $y \in \bb{Z}_N$ must remain occupied,
and have an arc of length at least $2$ immediately to its right,
for all times $t \leq T$. This contradicts the observation that
$Y^i$ passes through point of $\bb{Z}_N$ by time $T$. \qed

\textbf{Completion of proof of Proposition
\ref{prop:detBML}}\hspace{5pt} This is now immediate: if $m <
\frac{1}{2}N$ then, however the images $Y^1_t,Y^2_t,\ldots,Y^m_t$
are distributed in $\bb{Z}_N$, there will always be some empty arc
of length at least $2$, and so, by Lemma \ref{lem:nolongarcs}, the
system must attain speed one in finite time. \qed

Thus the system must self-organize from any (deterministic)
starting configuration with fewer than $\frac{1}{2}N$ cars.  A
related randomized question is:

\begin{qu}
Suppose we place a configuration of $m < \frac{1}{2}N$ cars on
$\bb{Z}_N \times \bb{Z}_N$ uniformly at random.  By Proposition
\ref{prop:detBML} the system will self organize to attain speed
one, but how many collisions will occur before it does so?
\end{qu}

The results above notwithstanding, problems can arise with as few
as $2N$ cars:

\begin{prop}\label{prop:stuck}
There is a configuration with $m$ cars which is stuck if and only
if $m \geq 2N$.
\end{prop}

\textbf{Proof}\hspace{5pt}
$\boldsymbol{(\Rightarrow)}$\hspace{5pt} Observe first that no
column can contain some red but no blue cars and no row can
contain some blue but no red cars, as then those cars would be
able to move freely. Now suppose the system is stuck. Then, in
particular, in every column there must be at least one blue car:
for there cannot be only red cars, and if there were no cars in
that column, then there would be blue cars in some column to the
left of it which are not blocked. Similarly, in every row there
must be at least one red car. Thus, there are at least $N$ red
cars and at least $N$ blue cars, so $m \geq 2N$.

$\boldsymbol{(\Leftarrow)}$\hspace{5pt} We need only witness a
configuration of $2N$ stuck cars: choose two adjacent SW-NE
diagonals, and occupy the lower entirely with red cars and the
upper entirely with blue cars. Now for any $m \geq 2N$, we can add
more cars that are also blocked to this configuration to deduce
the result. \qed

This leaves the following question for the deterministic problem:

\begin{qu}
Does the system necessarily self-organize to attain speed $1$ for
any $m \geq \frac{1}{2}N$?  Put differently, for which
$\frac{1}{2}N \leq m < 2N$ can the system stay forever below speed
$1$, even though (by Proposition \ref{prop:stuck}) it can never
get stuck?
\end{qu}

We remark that the idea underlying Proposition \ref{prop:stuck}
can also be brought to bear on a version of the model with a
(fairly sparse) random initial configuration:

\begin{prop}
If a configuration of $n\log n$ red cars and $n\log n$ blue cars
is selected in $\bb{Z}_N \times \bb{Z}_N$ uniformly at random,
then asymptotically almost surely the system never becomes stuck.
\end{prop}

\textbf{Proof}\hspace{5pt} We use again the observation that any
stuck configuration must contain at least one blue car in every
column and at least one red car in every row.  However, with $n
\log n$ cars of each colour, asymptotically almost surely there is
some column with no blue cars or some row with no red cars; this
gives the result. \qed

Finally, we see no harm in asking:

\begin{qu}
Could ideas such as those above be used to show that in the
original (random) BML model self organization occurs
(asymptotically almost surely) for $N^{1 +\alpha}$ cars for some
$0 < \alpha <1$?
\end{qu}

\end{document}